\documentclass[11pt]{amsart}

\title{Automorphism groups of relatively free groups  of   infinite rank:
 a survey}

\author{Vitali\u\i\ Roman'kov}
\address{Institute of Mathematics and Information Technologies\\Dostoevsky Omsk State  University}
\curraddr{}
\email{romankov48@mail.ru}

\usepackage{amsmath,amsthm,amsfonts,amssymb}
\usepackage{graphicx}
\usepackage{hyperref}
\usepackage[usenames]{color}

\theoremstyle{definition}

\newcounter{comcount}

\date{}

\begin{document}

\maketitle

\begin{abstract}
The paper is intended to be a survey on some topics within the framework of automorphisms of a relatively free groups of  infinite rank.   
We discuss such properties as {\it tameness,   primitivity, small index, Bergman property}, and so on.  
\end{abstract}


\begin{center}
\textit{Introduction}
\end{center}

Let $F_{\infty}$  be a free group of infinite rank, in particular, let $F_{\omega}$ be a free group of countably infinite rank.   In further of the paper $X_{\omega} =\{x_1, ..., x_i, ...\}$ be a basis of   $F_{\omega},$ and $X_n=\{x_1, ..., x_n\}$ be a basis of a free group $F_n$ of any finite rank $n.$  Thus $F_n$ is naturally embedded into $F_{\omega},$ and $F_n$ is naturally embedded into every group $F_m$ where $m \geq n.$ More generally, let $\Lambda $ be a set (finite or infinite) and $F_{\Lambda}$ be the free group of rank $|\Lambda |$ with  basis  $X_{\Lambda}=\{x_{\lambda}: \lambda \in \Lambda \}$. Then for each subset $\Xi $ of $\Lambda $ free group $F_{\Xi}$ is a subgroup of $F_{\Lambda},$  and for every $\Psi , \Xi \subseteq \Psi \subseteq \Lambda$, free group $F_{\Xi}$ is a subgroup of $F_{\Psi}$. 

For any variety of groups ${\mathcal C}$, let $V={\mathcal C}(F_{\Lambda })$  denote the verbal subgroup of $F_{\Lambda }$ corresponding
to ${\mathcal C}.$ (See \cite{Neu}  for information on varieties and related concepts.).  Then $G_{\Lambda } = F_{\Lambda }/V$ is  the free group of  rank $|\Lambda |$ 
in ${\mathcal C}$. In particular,  $G_{\omega }$ is the free group of countably infinite rank in ${\mathcal C}.$ Write $\bar{x}_i = x_iV$  for $i = 1,... ,  i, .. . \ .$ Then $X_V=\{\bar{x}_1,..., \bar{x}_i, ...\}$ is a basis of $G_{\omega}$. For each $\Lambda $ there is the standard homomorphism of $F_{\Lambda }$ onto $G_{\Lambda }$. Then for each subset $\Xi $ of $\Lambda $ free group $G_{\Xi}$ is a subgroup of $G_{\Lambda},$  and for every $\Psi , \Xi \subseteq \Psi \subseteq \Lambda$, free group $G_{\Xi}$ is a subgroup of $G_{\Psi}$.

If $\alpha $ is an automorphism of $G_{\Lambda }$ then $\{\alpha  (\bar{x}_{\lambda }) : \lambda \in \Lambda  \}$  is also a basis
of $G_{\Lambda }$ and every basis of $G_{\Lambda }$ has this form.

Any automorphism $\phi $ of $F_{ \Lambda }$ induces an automorphism $\bar{\phi} $  of $G_{\Lambda }$. Thus every basis of $F_{\Lambda }$ induces a basis of $G_{\Lambda }$. The converse
however is not always true; in general, there are automorphisms of $G_{\Lambda }$ which are not induced by automorphisms of $F_{\Lambda }$. See \cite{R1}, \cite{R2} for relevant results. 

An automorphism of $G_{\Lambda }$  which is induced by an automorphism of $F_{\Lambda }$ is called {\it tame}. 
If $\{g_{\xi} : \xi \in \Xi   \}$ are distinct elements 
of $G_{\Lambda }$ such that $\{g_{\xi} : \xi \in \Xi  \}$ is contained in a basis
of $G_{\Lambda }$ 
then $\{g_{\xi} : \xi \in \Xi \}$ is called a {\it primitive system} of $G_{\Lambda }$. 

Any primitive system $\{f_{\xi}: \xi \in \Xi \}$  of $F_{\Lambda }$   induces a primitive system of $G_{\Lambda }$ that is called {\it tame.} But, in general, not every primitive system of $G_{\Lambda }$ is induced by a primitive
system of $F_{\Lambda }$. We observe different tameness properties in the following Section 2.

An other topic of this paper is small index property. A countable first-order structure $M$ is said to have the small index property if every subgroup of the automorphism group Aut($M$) of index less than $2^{\aleph_0}$
contains the pointwise stabilizer C($U$) of a finite subset $U$ of the domain of $M $. In Section 3, we give results about small index property for 
relatively free groups of countably infinite rank.

Further in the paper, ${\mathcal A}$ denotes the variety of all abelian groups, ${\mathcal N}_k$ means the variety of all nilpotent groups of class $\leq k,$ and ${\mathcal A}^2$ stands for 
the variety  of all metabelian groups  (for any varieties ${\mathcal C}$  and ${\mathcal D}$, ${\mathcal C}{\mathcal D}$  denotes the variety of all groups with a normal subgroup in ${\mathcal C}$  and factor group in ${\mathcal D}$, thus ${\mathcal A}^2 = {\\mathcal A}{\mathcal A}$). We also denote by  $A_{\infty}, N_{k,\infty}$  and $M_{\omega}$ the free groups of the countably infinite rank in the varieties
${\mathcal A}, {\mathcal N}_k$ and ${\mathcal A}^2$, respectively. For any group $H$ and each positive integer $k$ we denote by $\gamma_k(H)$ the $k$th member of the low central series in $H.$ 
In particular, $\gamma_1(H)=H$ and $\gamma_2(H)= H'$, the derived subgroup of $H.$  

\begin{center}
\textit{Tameness}
\end{center}

In this section, we present known results about tame automorphisms of relatively free groups of  infinite rank. It is well known that every automorphism of $A_{\infty}$ can be lifted to an automorphism of $F_{\infty}$, thus is tame.
The following results also belong to this direction.

Theorem 1(Bryant and Macedonska \cite{BM}).
Let $F_{\infty}$  be a free group of infinite rank and let $V $ be a characteristic  
subgroup of $F_{\infty}$  such that $F_{\infty}/V$ is nilpotent. Then $G_{\infty}=F_{\infty}/V$  is relatively free group of infinite 
rank in a nilpotent variety and every automorphism of
$G_{\infty}$  is induced by an automorphism of $F_{\infty}$, thus is tame. 

If $F_{\infty}$  and $V$ are  as in the statement of the theorem then $V$ contains $\gamma_k(F_{\infty})$ for some
positive integer $k$.  Since V is characteristic in $F$, it follows by a result of 
Cohen \cite{Cohen1} that $V$ is fully characteristic in $F_{\infty}$. Thus $V$ defines 
nilpotent variety ${\mathcal N}$ of groups, then $V={\mathcal N}(F_{\infty})$ is the corresponding verbal
subgroup of $F_{\infty}$. Hence $F_{\infty}/V$ is a relatively free group in ${\mathcal N}$. 

To prove Theorem 1 the authors defined a property called
the {\it finitary lifting property} (see details below) and obtained the following  two results.

Proposition 1. Every nilpotent variety of groups has the finitary lifting
property.

Proposition 2. If ${\mathcal C}$  is any variety of groups with the finitary lifting
property and $F _{\infty}$ is a free group of infinite rank then every automorphism of
$F_{\infty}/{\mathcal C}(F_{\infty})$  is induced by an automorphism of $F_{\infty}$.

Let $F_{\infty}$  be a free group of infinite rank and let $\{x_{\lambda}: \lambda \in \Lambda\}$  be a basis of $F_{\infty}$. 
An automorphism $\phi$ of $F_{\infty}$  will be called {\it finitary} if
there is a finite subset  $U$ of this basis  such that  $\phi (x)=x$ 
for each free generator $x \not\in U$.
Let ${\mathcal C}$  be a variety of groups and write $V= {\mathcal C}(F_{\infty})$. Suppose that $\Gamma$ and $\Delta$
are subsets of $\Lambda$ such that $\Gamma \cap \Delta = \emptyset $, $\Delta $ is finite, and $\Lambda \setminus (\Gamma \cup \Delta)$ is infinite.
Let $\alpha $ be automorphism of $F_{\infty}/V$ such that 
$\alpha (x_{\lambda}V) = x_{\lambda}V$ for all $\lambda \in \Gamma .$ 
 We say that the triple $(\Gamma , \Delta , \alpha )$ can be lifted if there exists a finitary
automorphism $\xi $ of $F_{\infty}$ such that  $\xi (x_{\lambda}) = x_{\lambda}$ for all $\lambda$ in $\Gamma$  and $\xi (x_{\lambda}V) = \alpha (x_{\lambda}V)$
for all $\lambda \in \Delta $. Such a finitary automorphism $\xi$  is called a {\it lifting} of $(\Gamma , \Delta , \alpha ).$
We say that ${\mathcal C}$  has the {\it finitary lifting property} if, for every $F_{\infty}$ of infinite
rank, every triple $(\Gamma ,  \Delta , \alpha )$ can be lifted.

The theorem generalises some previously known results. The case where $V=\gamma_2(F_{\infty})$
is a result of Swan  (see \cite{Cohen2}).
A closely related result
had been obtained a few years earlier by Burns and Farouqi \cite{BF} who
proved that if $A_{\omega}(p)$ is a free abelian of exponent $p$ group of countably infinite rank and $p$ is a
prime number then every automorphism of $A_{\omega}(p)$ is induced by an
automorphism of $A_{\omega}$. In \cite{GM}, Gawron and Macedonska proved the discussed property in the cases $V = \gamma_i(F_{\omega})$ for $i = 3,4.$

For each positive integer $m$, we denote by ${\mathcal A}(m)$ the variety of all abelian groups of exponent dividing $m.$ Also we denote by ${\mathcal A}(0)$ the variety of all abelian groups ${\mathcal A}.$ 

Theorem 2 (Bryant and Groves \cite{BGr1}).
Let $m$ and $n$ be non-negative integers. Every automorphism a free group of infinite rank in the metabelian product variety ${\mathcal A}(m){\mathcal A}(n)$ is tame. In particular, every automorphism  of  $M_{\infty} $ is tame.

Theorem 3 (Bryant and Gupta \cite{BGup}). Let ${\mathcal C}$ be a variety  such that ${\mathcal A}^2 
\subseteq {\mathcal C} \subseteq {\mathcal N}_k{\mathcal A}$ for some $k,$ and $G_{\infty}$ be a free group of infinite rank in ${\mathcal C}.$ Then every automorphism of $G_{\infty}$ is tame.  

The following result generalizes Theorems 1, 2 and 3.

Theorem 4 (Bryant and Roman'kov \cite{BR1}). Let ${\mathcal C}$ be a subvariety of ${\mathcal N}_k{\mathcal A}$ for some $k.$ Let $G_{\infty}$ be a free group of infinite rank in ${\mathcal C}.$ Then every automorphism of $G_{\infty}$ is tame. 

The main ingredient in the proof of this theorem are the following  result that has its own interest.

Theorem  5 (Bryant and Roman'kov \cite{BR1}). Let ${\mathcal C}$ be a subvariety of ${\mathcal N}_k{\mathcal A}$, where $k \geq 1$. Let $n$ be a positive integer and write $l = 2^n(n+1) + 2k.$ Then every primitive system of $F_n/{\mathcal C}(F_n)$ is induced by some primitive system of $F_l.$ 

Above we presented some positive results about tameness of the automorphisms of relatively free groups of infinite rank. However, there are negative results for other varieties. 

Theorem 4 (Bryant and Groves \cite{BGr2}).
 Let ${\mathcal K} =$ var($K)$ be the variety generated by a non-abelian finite
simple group $K$, and $G_{\infty}$ is the free group of the countable infinite rank in ${\mathcal K}$. Then there is an automorphism of $G_{\infty}$ which is not induced by an automorphism of $F_{\infty}$.

\begin{center}
\textit{Small index property}
\end{center}

Hodges, Hodkinson, Lascar and Shelah  established in \cite{HHLS} that $\omega$-categorical and $\omega$-stable structures, and so called random graph have the small index property.  
 In \cite{BE}, Bryant and Evans   use the methods of the paper \cite{HHLS} to show that the free group of countably infinite rank and certain relatively free groups of countably infinite rank have the small index property.

Theorem 5 has some immediate consequences through the results of \cite{BE} and \cite{BR2}.

Theorem 6 (Bryant and Roman'kov \cite{BR1}). Let $F_{\omega }$ be a free group of countably infinite rank and let ${\mathcal C}$ be a subvariety of ${\mathcal N}_k{\mathcal A}$, where $k \geq 1$. Then $F_{\omega}/{\mathcal C}(F_{\omega}$ has the basis cofinality property and the small index property. The automorphism group Aut($F_{\omega}/{\mathcal C}(F_{\omega}$) is not the union of a countable chain of proper subgroups. Also, Aut($F_{\omega}/{\mathcal C}(F_{\omega}$) has no proper normal subgrou of index less than $2^{\aleph_0}$ and it is a perfect group. 

Recall that a group is called {\it perfect} if it equals its derived subgroup.

\begin{center}
\textit{Other  properties}
\end{center}

{\bf Completeness}. A group $G$ is said to be  complete if $G$ is centreless and every
automorphism of $G$ is inner. By the Burnside's criterion for a  centerless group $G$  its the automorphism group  Aut($G$) is complete if and
only if the subgroup Inn($G$) of all inner automorphisms of $G$ is a
characteristic subgroup of the group Aut($G$) (that is, preserved under the
action of all automorphisms of the group Aut($G$)). 

Theorem 7 (Tolstykh \cite{Tfree}. The automorphism group Aut($F_{\infty}$) of any  free group of infinite rank is complete.

This statement was derived from the following assertions:
\begin{itemize}
\item 
The family of all inner automorphisms of $F_{\infty }$  determined by powers of
primitive elements of $F_{\infty}$  is first-order definable in Aut($F_{\infty}$), hence 
Inn($F_{\infty}$) is a characteristic subgroup of Aut($F_{\infty}$).
\item The subgroup Inn($F_{\infty}$) is then first-order definable in Aut($F_{\infty}$). 
\end{itemize}

Theorem 8 (Tolstykh \cite{Tnil1}, \cite{Tnil2}). For any $k\geq 2$, the automorphism group Aut($N_{\infty , k}$) of any  free nilpotent group $N_{\infty , k}$ of infinite rank is complete. 

Note, that Inn(Aut($A_{\infty}$)) = Aut(Aut($A_{\infty}$)) (\cite{Tabcomm}, \cite{Tab}). Anyway Aut($A_{\infty}$) is not complete because it contains  a non-central the inverting authomorphism.  

In \cite{Tnil1}, this statement was proved for the case $k=2,$ and in \cite{Tnil2}, for the general case. 

Theorem 9 (Tolstykh \cite{Tr'}). 
Let $F_{\infty}$  be an infinitely generated free group, $R \leq  F_{\infty}'$ 
a fully characteristic subgroup of $F_{\infty }$ such that the quotient group $F_{\infty} /R$ is residually torsion-free nilpotent. Then the group Aut($F_{\infty}/R$)  is complete.

Corollary. Let $G_{\infty}$ be a free abelian-by-nilpotent (in particular metabelian or free solvable  of class $\geq 3$) relatively free group of infinite rank. Then the group Aut($G_{\infty}$)  is complete.

{\bf Generalized small index property}

Let $F$ be a relatively free algebra of infinite rank $\kappa $. We say that $F$ has the {\it generalized small index property} if any subgroup of   Aut($F$) of index at most $\kappa $  contains the pointwise stabilizer C($U$) of a subset $U$ of the domain of $F$ of cardinality less than $\kappa $. 

Theorem 10 (Tolstykh \cite{Tgsmall} Every infinitely generated free nilpotent  (in particular free abelian) group $N_{\infty}$ has the generalized small index property.

{\bf Bergman property}

A group $G$ is said to have the {\it Bergman property} (the property of {\it uniformity of finite width}) if given any generating $X$ with $X = X^{-1}$ of $G$, we have that $G = X^l$ for some natural $l$, that is, every element of $G$ is a product of at most $l$ elements of $X$. The property is named after Bergman, who found in \cite{Berg} that it is satisfied by all infinite symmetric groups. The first example of an infinite group with the Bergman property was apparently found by Shelah in the 1980s.

Theorem 10 (Tolstykh \cite{Tberg}).  The automorphism group Aut($F_{\omega}$) of the free group $F_{\omega}$ of countably infinite rank has the Bergman property.

Theorem 11 (Tolstykh \cite{Tberg}). For any positive integer $k$, the automorphism group Aut($N_{\infty , k}$) of any  free nilpotent group $N_{\infty , k}$  of infinite rank has the Bergman property.

Some other discussion on the automorphism groups of free relatively free groups can be found in survey \cite{Tnov}.

\end{document}